\documentclass{amsart}

\usepackage{amssymb, amsthm}

\newtheorem{thm}{Theorem}[section]
\newtheorem{lem}[thm]{Lemma}
\newtheorem{cor}[thm]{Corollary}

\theoremstyle{definition}

\newtheorem{ex}[thm]{Example}

\theoremstyle{remark}
\newtheorem{rem}[thm]{Remark}

\begin{document}

\title{Failure of F-purity and F-regularity in certain rings of invariants}

\author{Anurag K. Singh}

\address{Department Of Mathematics, The University Of Michigan, East Hall, 525
East \linebreak University Ave., Ann Arbor, MI 48109-1109}

\email{binku@math.lsa.umich.edu}


\maketitle

\section{Introduction}

Let ${\mathbb F}_q$ be a finite field of characteristic $p$, $K$ a field 
containing it, and $R=K[X_1, \dots, X_n]$ a polynomial ring in $n$ variables. 
The general linear group $GL_n({\mathbb F}_q)$ has natural action on $R$ by 
degree preserving ring automorphisms. L.~E.~Dickson showed that the subring of 
elements which are fixed by this group action is a polynomial ring, 
\cite{Dickson}, though for an arbitrary subgroup $G$ of $GL_n({\mathbb F}_q)$, 
the structure of the ring of invariants $R^G$ 
may be rather mysterious. If the order of the group $|G|$ is relatively prime
to the characteristic $p$ of the field, there is an $R^G$--linear retraction 
$\rho :R \to R^G$, the {\it Reynolds operator}\/. This retraction makes $R^G$ a 
direct
summand of $R$ as an $R^G$--module, and so $R^G$ is F--regular.  However when 
the characteristic $p$ divides $|G|$, this method no longer applies, and the 
ring of invariants $R^G$ need not even be Cohen--Macaulay. M.--J.~Bertin showed 
that when $R$ is a polynomial ring in four variables and $G$ is the cyclic group
with four elements which acts by permuting the variables in cyclic order, then 
the ring of invariants $R^G$ is a unique factorization domain which is not 
Cohen--Macaulay, providing the first example of such a ring, \cite{Bertin}. 
Related work and bounds on the depth of $R^G$ can be found in the work of 
R.~M.~Fossum and P.~A.~Griffith, see \cite{FossumGriffith}. More recently
D.~Glassbrenner studied the invariant subrings of the action of the alternating
group $A_n$ on a polynomial ring in $n$ variables over a field of characteristic
$p$, constructing examples of F--pure rings which are not F--regular,
\cite{Glthesis, Gl}. Both these families of examples study rings of 
invariants of
$K[X_1, \dots, X_n]$ under the action of a subgroup $G$ of the symmetric group
on $n$ elements, i.e., an action which permutes the variables, and Glassbrenner 
shows that for such a group the ring of invariants is F--pure, see 
\cite[Proposition 0.6.7]{Glthesis}. 

\medskip

We shall construct examples which demonstrate that the ring of invariants 
for the natural action of a subgroup $G$ of $GL_n({\mathbb F}_q)$ need not be 
F--pure. We shall obtain such examples with the group $G$ being the symplectic 
group over a finite field. These non F--pure invariant subrings are always
complete intersections, and are actually hypersurfaces in the case of  
$G=Sp_4({\mathbb F}_q) < GL_4({\mathbb F}_q)$ acting on the polynomial ring
$R=K[X_1,X_2,X_3,X_4]$. These examples are particularly interesting if one is
attempting to interpret the Frobenius closures and tight closures of ideals as
contractions from certain extension rings, since we have an ideal generated by
a system of parameters and the socle element modulo this ideal is being forced 
into the expansion of the ideal to a module--finite extension ring which is a 
separable, in fact Galois, extension. This element is also forced into an 
expanded ideal in a linearly disjoint purely inseparable extension, being in 
the Frobenius closure of the ideal. It is noteworthy that the element can be 
forced into expanded ideals in two such different ways.

\medskip

Our results depend on the work of D.~Carlisle and P.~Kropholler where they show 
that the ring of invariants under the natural action of the symplectic group on 
a polynomial ring is a complete intersection, \cite{CK}. We obtain the precise 
equations defining these complete intersections in some examples using the
program Macaulay, and in some other cases collect enough information to display 
that the invariant subrings are not F--pure. 

\medskip

The second part of this paper deals with the alternating group $A_n$ acting on 
the polynomial ring $R=K[X_1, \dots, X_n]$ by permuting the variables. We shall 
assume that the characteristic $p$ of $K$ is an odd prime, and denote by 
$R^{A_n}$, the invariant subring of this action. Since $R^{A_2}$ is a polynomial
ring we shall always assume $n \ge 3$. If the order 
of the group $|A_n| = \frac{1}{2}(n!)$ is relatively prime to the 
characteristic $p$ of the field, the Reynolds operator makes $R^{A_n}$ a direct 
summand of $R$ as an $R^{A_n}$--module, and in the language of tight closure, 
the existence of such a retraction is equivalent to the ring $R^{A_n}$ 
being F--regular, see Lemma \ref{L:checksocle}. When $p$ divides $n$ or $n-1$,
Glassbrenner has shown that the invariant subring $R^{A_n}$ is no longer 
F--regular, see \cite[Proposition 1.2.5]{Glthesis}. We shall extend this result 
by showing that $R^{A_n}$ is F--regular if and only if $p$ does not divide $|A_n|$. 
 
\medskip

The author wishes to thank Melvin Hochster for several interesting discussions. 

\section{F--purity and F--regularity}

We recall some basic notation and definitions from 
\cite{HHjams, HHbasec, HHjalg}. 

\medskip

Let $R$ be a Noetherian ring of characteristic $p > 0$. We shall always use
the letter $e$ to denote a variable nonnegative integer, and $q$ to denote the
$e\,$th power of $p$, i.e., $q=p^e$. For an ideal 
$I=(x_1, \dots, x_n) \subseteq R$ we let $I^{[q]}=(x_1^q, \dots, x_n^q)$. 

\medskip

For an element $x$ of $R$, we say that $x \in I^F$, the 
{\it Frobenius closure }\/ of $I$, if there exists $q=p^e$ such that
$x^q\in I^{[q]}$. We shall say that the ring $R$ is 
{\it F--pure }\/ if for all ideals $I$ of $R$ we have  $I^F = I$.

\medskip

We shall denote by $R^{\text{o}}$ the complement of the union of the minimal primes of
$R$. For an ideal $I \subseteq R$ and an element $x$ of $R$, we say that 
$x \in I^*$, the {\it tight closure }\/ of $I$, if there exists 
$c \in R^{\text{o}}$ such that $cx^q \in I^{[q]}$ for all $q=p^e \gg 0$. If 
$I=I^*$ for all ideals $I$ of $R$, we say $R$ is {\it weakly F--regular}. $R$
is called {\it F--regular}\/ if every localization is weakly F--regular. These
two notions are known to be the same if $R$ is Gorenstein, 
\cite[Corollary 4.7]{HHbasec}.

\section{Symplectic invariants}

We shall summarize in this section the results of Carlisle and Kropholler as 
presented in \cite{Benson}. Let ${\mathbb F}_q$ be a finite field of 
characteristic $p$, and $K$ an infinite field containing it. L.~E.~Dickson 
showed that the ring of invariant forms
under the natural action of $GL_n({\mathbb F}_q)$ on the polynomial ring 
$R=K[X_1, \dots, X_n]$ is a graded polynomial algebra on the algebraically 
independent generators $c_{n,i}$, where the $c_{n,i}$ are the coefficients in 
the equation
$$
\prod_{v \in {\mathbb F}_q[X_1, \dots, X_n]}(T-v) = T^{q^n}-c_{n,n-1}T^{q^{n-1}}+c_{n,n-2}T^{q^{n-2}}- \dots
+(-1)^nc_{n,0}T.
$$
When working with a fixed polynomial ring $R=K[X_1, \dots, X_n]$, we shall drop 
the first index, and write the generators of $R^{GL_n({\mathbb F}_q)}$ as 
$c_0, \dots, c_{n-1}$, the {\it Dickson invariants}.
It is clear that for any subgroup $G$ of $GL_n({\mathbb F}_q)$, the ring of 
invariants $R^G$ is a module--finite extension of the polynomial ring
$R^{GL_n({\mathbb F}_q)} = K[c_0, \dots, c_{n-1}]$. 

\medskip

Let $V$ be a vector space of dimension $2n$ over the field ${\mathbb F}_q$, on 
which we have a non--degenerate alternating bilinear form $B$. We may choose a 
basis $e_1, \dots, e_{2n}$ for $V$, such that $B$ is given 
by
$$
B(\sum a_i e_i, \sum b_j e_j)
= a_1b_2 - a_2b_1 + \dots + a_{2n-1}b_{2n} - a_{2n}b_{2n-1}.
$$
The symplectic group $G=Sp_{2n}({\mathbb F}_q)$ is the subgroup of 
$GL_{2n}({\mathbb F}_q)$ consisting of the elements which preserve $B$. We 
consider the natural action of $G$ on $R = K[X_1, \dots, X_{2n}]$. In addition 
to the Dickson invariants, it is easily seen that $R^G$ must contain
$$
\xi_i = X_1X_2^{q^i}-X_2X_1^{q^i}+ \dots 
+ X_{2n-1}X_{2n}^{q^i}-X_{2n}X_{2n-1}^{q^i}.
$$
Carlisle and Kropholler show that the Dickson invariants $c_0, \dots, c_{2n-1}$
along with the above $\xi_1, \dots, \xi_{2n}$ form a generating set for $R^G$, 
and that there are $2n$ relations, i.e., that $R^G$ is a complete intersection. 
One may eliminate $c_0, \dots, c_{n-1}$ and $\xi_{2n}$ using $n+1$ of these 
relations, after which the remaining $n-1$ relations are 
$$
\sum_{j=0}^{i-1} (-1)^j \xi_{i-j}^{q^j}c_j 
= \sum_{j=i+1}^{2n} (-1)^j \xi_{j-i}^{q^i}c_j 
$$
where $1 \le i \le n-1$ and $c_{2n}=1$. Their results furthermore show that
$c_0 \in K[\xi_1, \dots, \xi_{2n-1}]$ which is, in fact, a polynomial ring.

\section{Rings of invariants which are not F--pure}

We shall first show that the ring of invariants of $G=Sp_4({\mathbb F}_q)$ 
acting on the polynomial ring $R=K[X_1,X_2,X_3,X_4]$ is not F--pure when $q=2$ 
or $3$. Note that $Sp_2({\mathbb F}_q)$ is the same as $SL_2({\mathbb F}_q)$, 
and so the ring of invariants in that case is a polynomial ring. 

\begin{ex}
Let $R=K[X_1,X_2,X_3,X_4]$ and $G=Sp_4({\mathbb F}_q)$ be the symplectic group 
with its natural action on $R$. In the notation of the previous section, 
$R^G = K[c_2, c_3, \xi_1, \xi_2, \xi_3]$, where the only relation is
$$
\xi_1 c_0 = \xi_1^q c_2 - \xi_2^q c_3 + \xi_3^q.
$$
We need to determine $c_0$ as an element of $K[\xi_1, \xi_2, \xi_3]$. 
When $q=2$, it can be verified that $c_0 = \xi_1^5 + \xi_2^3 + \xi_3 \xi_1^2$,
and so 
$$
\xi_3^2 = \xi_1^6 + \xi_1 \xi_2^3 + \xi_1^3 \xi_3 + \xi_1^2 c_2 + \xi_2^2 c_3,
$$
by which $\xi_3 \in ((\xi_1, \xi_2)R^G)^F$. Since 
$\xi_3 \notin (\xi_1, \xi_2)R^G$, the ring $R^G$ is not F--pure.

\medskip

In the case $q=3$, $c_0$ can be expressed as an element of 
$K[\xi_1, \xi_2, \xi_3]$ by the equation
$$
c_0 = \xi_2^8 + \xi_3 \xi_1^3 \xi_2^4 + \xi_1^6 \xi_3^2 + \xi_1^{10} \xi_2^4
-\xi_1^{13} \xi_3 + \xi_1^{20}.
$$ Once again we see that 
$\xi_3 \in ((\xi_1, \xi_2)R^G)^F$, and so $R^G$ is not F--pure. 

\medskip

Computations with Macaulay helped us determine the precise equations in these
examples. 
\label{nfpinvar}
\end{ex}

\begin{thm}
Let ${\mathbb F}_q$ be a finite field of characteristic $p$, and $K$ an
infinite field
containing it. Let $G=Sp_{2n}({\mathbb F}_q)$ be the symplectic group with its
natural action on the polynomial ring $R = K[X_1, \dots, X_{2n}]$. If $n \ge 2$
and $q \ge 4n-4$, then the ring of invariants $R^G$ is not F--pure. 
\label{T:nfpinvar}
\end{thm}

\begin{proof}
In the notation of the previous section, the ring of invariants is 
$R^G=K[c_n, \dots, c_{2n-1},\xi_1, \dots, \xi_{2n-1}]$, where there are exactly
$n-1$ relations, as stated before. Using the relation with $i=1$, we see that
$$
\xi_{2n-1}^q \in (\xi_1^q, \dots, \xi_{2n-2}^q, \xi_1c_0)R^G,
$$
whereas $\xi_{2n-1} \notin (\xi_1, \dots, \xi_{2n-2})R^G$. 

\medskip

If $R^G$ is indeed
F--pure, $\xi_{2n-1}^q \notin (\xi_1^q, \dots, \xi_{2n-2}^q)R^G$, and so the 
expression of $c_0$ as an element of $K[\xi_1, \dots, \xi_{2n-1}]$ 
must have a monomial of the form 
$\xi_1^{a_1} \xi_2^{a_2} \dotsm \xi_{2n-1}^{a_{2n-1}}$,
with $a_1 \le q-2$ and $a_2, \dots, a_{2n-1} \le q-1$. Equating degrees, we
have
\begin{align*}
\deg c_0 = q^{2n}-1 &= a_1(q+1)+a_2(q^2+1)+ \dots + a_{2n-1}(q^{2n-1}+1) \\
        &= \sum_{i=1}^{2n-1}a_i + \sum_{i=1}^{2n-1}a_iq^i.
\end{align*}

Examining this modulo $q$, we get that $\sum_{i=1}^{2n-1}a_i = \lambda q -1$,
where the bounds on $a_i$ show that $1 \le \lambda \le 2n-2 < q$.
Substituting this, we get $q^{2n} = \lambda q + \sum_{i=1}^{2n-1}a_iq^i$.
Working modulo $q^2$, we see that $a_1 = q - \lambda$, and continuing this way
we get that $a_2, \dots, a_{2n-1} = q-1$. Hence
$$
q^{2n}-1 = (q-\lambda)(q+1)+(q-1)(q^2+1)+ \dots + (q-1)(q^{2n-1}+1),
$$
which simplifies to give $\lambda(q+1) = 2nq-2n-q+3$. Since $\lambda \le 2n-2$,
this implies that $q \le 4n-5$, a contradiction.

\medskip

Hence $R^G$ is not F--pure. In particular 
$\xi_{2n-1} \in ((\xi_1, \dots, \xi_{2n-2})R^G)^F$, the Frobenius closure. 
\end{proof}

\begin{cor} 
The ring of invariants of the symplectic group 
$G=Sp_4({\mathbb F}_q)$ acting on the polynomial ring $R=K[X_1,X_2,X_3,X_4]$ is not 
F--pure.
\label{C:nfpinvar}
\end{cor} 

\begin{proof}
We have, in the examples above, treated the case where $q=2$ or $3$. When 
$q \ge 4$, the result follows from the previous theorem. 
\end{proof}

\section{Rings of invariants of the alternating group}

The invariant subring under the natural action of the alternating group $A_n$ is 
$R^{A_n} = K[e_1, \dots, e_n, \Delta]$ where $e_i$ is the elementary symmetric 
function of degree $i$ in $X_1, \dots, X_n$, 
and $\Delta = \prod_{i > j}(X_i-X_j)$. The element $\Delta$ is easily seen to be
fixed by all even permutations of $X_1, \dots, X_n$, though not by odd
permutations. However its square, $\Delta^2$, is fixed by all permutations, and
so is a polynomial in the algebraically independent elements $e_1, \dots, e_n$. 
Consequently the invariant subring $R^{A_n}$ is a hypersurface, in particular
it is Gorenstein. The elements $e_1, \dots, e_n$ are an obvious choice as a 
homogeneous system of parameters for $R^{A_n}$, and  the one--dimensional socle
modulo this system of parameters is generated by $\Delta$. 

\begin{lem} 
With the above notation, the following are equivalent:

\item $(1)$\quad $R^{A_n} = K[e_1, \dots, e_n, \Delta]$ is F--regular. 

\item $(2)$\quad $R^{A_n}$ is a direct summand of $R=K[X_1, \dots, X_n]$ as an
$R^{A_n}$--module. 

\item $(3)$\quad $\Delta \notin (e_1, \dots, e_n)R$.
\label{L:checksocle}
\end{lem}

\begin{proof} 
$(1)\Rightarrow (2)$\quad By \cite[Theorem 5.25]{HHjalg}, an F--regular ring
is a direct summand of any module--finite extension ring.

$(2)\Rightarrow (3)$\quad Since $R^{A_n}$ is a direct summand of $R$, we have
$$
(e_1, \dots, e_n)R \cap R^{A_n} = (e_1, \dots, e_n)R^{A_n}.
$$ 

$(3)\Rightarrow (1)$\quad The elements $e_1, \dots, e_n$ form a system of parameters 
for the Gorenstein ring $R^{A_n}$ and $\Delta$ is the socle generator modulo
this system of parameters. If $\Delta$ is in the tight closure of
$(e_1, \dots, e_n)R^{A_n}$, then 
$\Delta \in (e_1, \dots, e_n)R^*=(e_1, \dots, e_n)R$. Hence $\Delta$ cannot
be in the tight closure of $(e_1, \dots, e_n)R^{A_n}$, by which $R^{A_n}$ is
F--regular.
\end{proof} 

Consequently our aim is to establish that $\Delta \in (e_1, \dots, e_n)R$,
whenever $p$ divides $|A_n|$. We shall henceforth denote this ideal by
$I = (e_1, \dots, e_n)R$.

\begin{lem}
Let $T_j^i$ denote the sum of all monomials of degree $i$ in the variables
$X_j, \dots, X_n$. Then $T_j^i \in I$ whenever $i \ge j \ge 1$. In particular, 
$T_i^i \in I$ for all $1 \le i \le n$. 
\label{L:alt1}
\end{lem}

\begin{proof}
Observe that $T_j^i =T_{j-1}^i - X_{j-1}T_{j-1}^{i-1}$. Given $T_j^i$ with 
$i \ge j \ge 1$, we may use this formula to rewrite $T_j^i$ as a sum of terms
which are multiples of $T_1^i$. Since $T_1^i$ is the sum of all the monomials 
of degree $i$ in $X_1, \dots, X_n$, it is certainly an element of $I$, and so
$T_j^i \in I$. 
\end{proof}

\begin{lem}
The ideal $I = (e_1, \dots, e_n)R$ generated by the elementary symmetric
functions contains the elements: $X_n^n$, \ $X_n^{n-1}X_{n-1}^{n-1}$, \  
$X_n^{n-1}X_{n-1}^{n-2}X_{n-2}^{n-2}$, \ $\dots$, \ 
$X_n^{n-1}X_{n-1}^{n-2} \dotsm X_i^{i-1}X_{i-1}^{i-1}$, \ $\dots$, \
$X_n^{n-1}X_{n-1}^{n-2} \dotsm X_2X_1$.
\label{L:alt2}
\end{lem}

\begin{proof}
We shall use the fact that $T_i^i \in I$ for $1 \le i \le n$, 
Lemma \ref{L:alt1}. This already
says that $X_n^n = T_n^n \in I$, and since $I$ is symmetric in the $X_i$,
we also have $X_{n-1}^n \in I$. Next, $X_{n-1}^{n-1}T_{n-1}^{n-1} \in I$,
but examining this using $X_{n-1}^n \in I$ we see that 
$X_n^{n-1}X_{n-1}^{n-1} \in I$. We proceed by induction.

\medskip

Since $T_{i-1}^{i-1} \in I$, we know that
$X_n^{n-1}X_{n-1}^{n-2} \dotsm X_i^{i-1} T_{i-1}^{i-1} \in I$, 
but using the inductive hypothesis this gives 
$$
X_n^{n-1}X_{n-1}^{n-2} \dotsm X_i^{i-1}X_{i-1}^{i-1} \in I.
$$
\end{proof}

\begin{lem}
In the above notation, 
$\Delta \equiv (n!) X_n^{n-1}X_{n-1}^{n-2} \dotsm X_2 \pmod I$.
\label{L:alt3}
\end{lem}

\begin{proof}
Let $\delta_r = (X_r-X_1)(X_r-X_2) \dotsm (X_r-X_{r-1})$. Then
$\Delta = \delta_n \delta_{n-1} \dotsm \delta_2$. We shall show that
$\delta_r \equiv r X_r^{r-1} \pmod {I + (X_{r+1},\dots, X_n)R}$ for 
$2\le r \le n$.
Note that for $r=n$, this says $\delta_n \equiv n X_n^{n-1} \pmod I$.

\medskip

Fix $r$, where $2\le r \le n$. Let $f_i$ be the elementary symmetric 
function of degree $i$ in the variables $X_1, \dots, X_{r-1}$. Then 
$$
f_i \equiv (-X_r)f_{i-1} \pmod {I + (X_{r+1},\dots, X_n)}R,
$$ 
and using this repeatedly, we see 
$$
f_i \equiv (-X_r)^i \pmod J \quad \text{ where }\quad J =I+ (X_{r+1},\dots, X_n)R.
$$
Consequently 
\begin{align*}
\delta_r &= (X_r-X_1)(X_r-X_2) \dotsm (X_r-X_{r-1}) \\
&=
X_r^{r-1} - X_r^{r-2}(X_1+\dots +X_{r-1}) + \dots + 
(-1)^{r-1}X_1\dotsm X_{r-1} \\
&\equiv 
X_r^{r-1} - X_r^{r-2}f_1 + \dots + (-1)^{r-1}f_{r-1} 
\pmod J \\
&\equiv 
X_r^{r-1} - X_r^{r-2}(-X_r) + \dots + (-1)^{r-1}(-X_r)^{r-1} 
\pmod J \\
&\equiv r X_r^{r-1} \pmod J. 
\end{align*}

Since $X_n^n \in I$, when evaluating the term $\delta_n \delta_{n-1} \pmod I$, 
it is enough to consider $\delta_{n-1} \pmod {I + X_nR}$, and using this we get 
$$
\delta_n \delta_{n-1} \equiv n(n-1)X_n^{n-1}X_{n-1}^{n-2} \pmod I.
$$
Proceeding in this manner, one obtains from the above calculations that 
$$
\Delta = \delta_n \delta_{n-1}\dotsm \delta_2 
          \equiv (n!)X_n^{n-1}X_{n-1}^{n-2} \dotsm X_2  \pmod I.
	  $$
The point is that since 
$$
\delta_n \delta_{n-1}\dotsm \delta_r 
        \equiv n(n-1)\dotsm(r) X_n^{n-1}X_{n-1}^{n-2} \dotsm X_r^{r-1} \pmod I,
$$	
we have $\delta_n \delta_{n-1}\dotsm \delta_r(X_r, \dots, X_n) \subseteq I$,
by Lemma \ref{L:alt2} and so when evaluating the product 
$\delta_n \delta_{n-1}\dotsm \delta_{r-1} \pmod I$, 
one need only consider the element $\delta_{r-1}$ modulo the ideal
$I +(X_r,\dots, X_n)R$. 
\end{proof}

We are now ready to prove the main result of this section.

\begin{thm}
Let $R=K[X_1, \dots, X_n]$ be a polynomial ring in $n$ variables over a field 
$k$ of characteristic $p$, an odd prime, and let the alternating group $A_n$ act
on $R$ by permuting the variables. Then the invariant subring $R^{A_n}$ is 
F--regular
(equivalently, $R^{A_n}$ is a direct summand of $R$) if and only if the order 
of the group $|A_n| = \frac{1}{2}(n!)$ is relatively prime to $p$.
\end{thm}

\begin{proof}
As we noted, it suffices to show
that $\Delta \in I = (e_1, \dots, e_n)R$. By Lemma \ref{L:alt3}, 
$\Delta \equiv (n!) X_n^{n-1}X_{n-1}^{n-2} \dotsm X_2 \pmod I$, and so the 
result follows. 
\end{proof}

\begin{rem}
\cite[Proposition 0.6.7]{Glthesis} shows that $R^{A_n}$ is always F--pure. 
Consequently when the characteristic $p$ of the field $K$ is an odd prime 
dividing $|A_n|$, $R^{A_n}$ is an F--pure ring which is not F--regular.
\end{rem}

\end{document}